\documentclass[11pt, a4paper]{article}
\usepackage{adjustbox}
\usepackage{srcltx}
\usepackage{amsmath,latexsym,amsfonts,amssymb}
\usepackage{pstricks}
\usepackage[latin1]{inputenc}
\usepackage{epsfig}
\usepackage{fancyhdr}
\usepackage{graphics,graphicx}
\usepackage{xspace}
\usepackage{enumerate}
\usepackage{rotating}
\usepackage{multirow}
\usepackage[linesnumbered,ruled,vlined,,procnumbered]{algorithm2e}
\usepackage{float}
\usepackage{pifont}
\usepackage{tikz}
\usepackage{pgfplots}
\usepackage[round]{natbib}
\usepackage{color}
\usepackage{tablefootnote}
\usepackage{authblk}
\usepackage{amssymb}
\usepackage{url}

\usepackage{subcaption}
\usepackage{makecell}
\renewcommand\theadfont{\normalsize\scshape}

\newcommand{\dCal}{\mathcal{D}}

\newcommand{\cCal}{\mathcal{C}}
\newcommand{\pCal}{\mathcal{P}}
\newcommand{\lCal}{\mathcal{L}}

\newcommand{\rCal}{\mathcal{R}}
\newcommand{\iCal}{\mathcal{I}}

\newcommand{\baralpha}{\overline{\alpha}}

\makeatletter
\AtBeginEnvironment{function}{%
  \let\c@algocf\c@function
}
\makeatother

\begin{document}
    
    \title{
        On Finding the Community with Maximum Persistence Probability}
    
    \author[(1)]{Alessandro Avellone}
    \author[(2)]{Stefano Benati}
    \author[(1)]{Rosanna Grassi}
    \author[(1)]{Giorgio Rizzini}
    
    \affil[(1)]{\small{Dipartimento di Statistica e Metodi Quantitativi, Universit\`{a} degli Studi di Milano - Bicocca, Milano}}
    
    \affil[(2)]{\small{Dipartimento di Sociologia e Ricerca Sociale, Universit\`{a} degli Studi di Trento, Italy, stefano.benati@unitn.it}}
    
    
    
    \date{}

    \maketitle
    
\begin{abstract}
The persistence probability is a statistical index that has been proposed to detect one or more communities embedded in a network. Even though its definition is straightforward, e.g, the probability that a random walker remains in a group of nodes, it has been seldom applied possibly for the difficulty of developing an efficient algorithm to calculate it.
Here, we propose a new mathematical programming model to find the community with the largest persistence probability. The model is integer fractional programming, but it can be reduced to mixed-integer linear programming with an appropriate variable substitution. Nevertheless, the problem can be solved in a reasonable time for networks of small size only, therefore we developed some heuristic procedures to approximate the optimal solution. First, we elaborated 
a randomized greedy-ascent method, taking advantage of a peculiar data structure to generate feasible solutions fast. After analyzing the greedy output and determining where the optimal solution is eventually located, we implemented improving procedures based on a local exchange, but applying different long term diversification principles, that are based on variable neighborhood search and random restart. Next, we applied the algorithms on simulated graphs that reproduce accurately the clustering characteristics found in real networks to determine the reliability and the effectiveness of our methodology. Finally, we applied our method to two real networks, comparing our findings to what found by two well-known alternative community detection procedures.
\end{abstract}
    

\section{Introduction}

The analysis of real networks, as they emerged as a structural model in disciplines as different as biology, economics, social sciences, engineering and so on, brought about a growing and thriving interest in developing new tools and methods to uncover the networks hidden characteristics, such as their communities, see \cite{Fortunato2016}, their core-periphery structure, see \cite{tang2019},  their node centrality, see \cite{das2018}. It is beyond our possibility to list all the contributions that mathematical programming provided to the field, as many networks features to find can be expressed by optimization models. Here we focus on a specific measure, the \textit{persistence probability}, that has been proposed to find one or more communities embedded in a network, see \cite{piccardi2011}. Loosely speaking, given a subset of nodes, its persistence probability is the probability that a random walker, located by chance in one of these nodes and moving randomly across the links, will remain in another node of the subset. 
The ratio is that this statistic should be able to detect nodes as well connected with each other to form a community: the highest the persistence, the highest fraction of links are directed towards internal nodes to the detriment of the external ones. The measure aroused some interest among scholars: for example, it has been used in \cite{della2013} to detect the core-periphery structures in many real networks such as the Karate Club, the co-authorship, the proteins and the World Trade networks. Further analyses of the World Trade through persistence are available in \cite{Piccardi2012} and they have been used to identifies the \textit{locali} (the local mobs) of the \textit{n'drangheta} criminal networks in \cite{Calderoni2017}.

The persistence probability has a clear and appealing definition and it is flexible enough to be applied for both community detection and core-periphery analysis. However, its application is still undervalued, maybe due to the fact that computational methods have not been developed yet, at least to the best of our knowledge. In this contribution,  we propose a new mathematical programming model to find the community with the largest persistence. One of the difficulty of the model is imposing connectivity on the node subset, but it is resolved using ideas developed in \cite{Benati2017}.  Next, the integer fractional objective function is reduced to mixed-integer linear programming using appropriate variable substitutions. Unfortunately, but predictability as the problem is NP-hard, the problem can be solved exactly only when the network size is small, therefore we developed some heuristic procedures to approximate the optimal solution. Peculiar of the maximum persistence problem is that the correct subset/community size $k$ must be known in advance, a case that rarely occurs in practice. Indeed, the correct value of $k$ can be guessed only after the analysis of various heuristic outcomes. Therefore, the implementation of a heuristic must consider that it should be fast enough to run the preliminary analysis to fix $k$ in a reasonable time. To the purpose, we elaborated on a randomized greedy-ascent method, proposed previously for a similar problem in \cite{Benati2022}. That algorithm takes advantage of a peculiar data structure coming from clustering problems, that generate feasible solutions fast. After determining the right value of $k$, the greedy outcome can be improved by local exchange and long-term diversification strategies. Here, we adopt diversification based on variable neighborhood search and random restart, but with some variation due to the problem structure. Variable neighborhood search has been implemented of the reduction of the cluster to its spanning tree, random restart is controlled by preliminary diversification.

We test the whole methodology, its accuracy and computational time, on graphs simulated through the procedure proposed in \cite{lancichinetti2008}, as was done also in \cite{piccardi2011}. This procedure simulates synthetic networks with the same characteristics found in real networks, therefore they are a severe and realistic benchmark. As it can be seen, the right size $k$ is often correctly determined after the greedy, and then the diversification heuristic improves the incumbent solution (when possible) in a fast way and hidden communities embedded in the network are detected. In the end, we apply our method to two real networks to test its ability in identifying communities, comparing its results with what found by two alternative methods, e.g., the Walktrap and the Louvain, see \cite{Blondel2008,pons2006}, and we will see that the use of the persistence complements well the findings of the other methods.

The paper is organized as follows. In Section \ref{Persistence_prob} the definition of persistence probability is formalized. In Section \ref{problemformulation}, finding the node subset with maximum persistent is formulated as an optimization problem, that after some modification is turn into mixed-integer linear programming. In Section \ref{Heuristic}, some heuristic algorithms are introduced: the first is a greedy procedure with some randomized steps, able to calculate the optimal persistence for subsets of varying size $k$, next the greedy results are improved with the interchange heuristic and some version of long-term diversification. In Section \ref{test}, computational tests are carried on to explain how to use the persistence probability and what are the best algorithms to find it. Finally, in Section \ref{real_analysis}, some empirical experiments on two real networks are explained. Conclusions follow (Section \ref{conclusions}).

\section{The Persistence probability}
\label{Persistence_prob}
\begin{figure}[htp]
	    	\centering
		\includegraphics[scale=0.5]{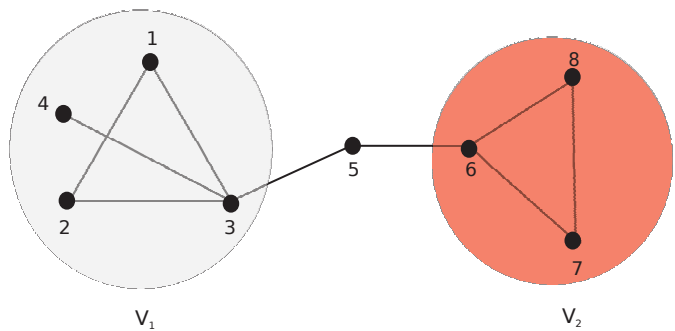}
		\caption{Different persistence probabilities of the subsets $V_1$ and $V_2$ for a graph $G$. The persistence probabilities are $\alpha(V_1)=\frac{4}{5}$ and $\alpha(V_2)=\frac{3}{4}$, respectively.
		}
			\label{fig:eesempiosemplice}
	\end{figure}

Let $G=(V,E)$  be a simple, undirected and connected graph (or network) where $V$ is the set of nodes and $E$ is the set of the edges (or links). Let $n = |V|$ be the cardinality of $V$. The square binary matrix $\mathbf{A}$ of elements $a_{ij}$ $i,j=1,...,n$ represents the adjacency relationships between nodes. Being the network undirected, $\mathbf{A}=\mathbf{A}^T$, where $\mathbf{A}^T$ is the transpose of $\mathbf{A}$. Consider a node subset $V_C \subseteq V$ and assume that the subgraph $G_C$, induced \footnote{The induced subgraph $G_C$ is the graph whose vertex set is $V_C$ and whose edge set consists of all the edges in $E$ that have both endpoints in $G_C$.} by $V_C$ is connected. Let $E_C$ be the edge set of the subgraph.
In \cite{piccardi2011} and \cite{della2013}, the persistence probability $\alpha(V_C)$ is proposed  
as a measure of cohesiveness of  subset $V_C$. Formally, $\alpha(V_C)$ is defined as:

\begin{equation}
\alpha(V_C) = \frac{\sum_{(i,j) \in E_C} a_{ij}}{\sum_{i \in V_C}\sum_{j \in V} a_{ij}},
\label{formula_alfa}
\end{equation}
    
\noindent
expressing the ratio between the number of links connecting nodes inside $V_C$, e.g. the internal links,
and all the links emanating from $V_C$, e.g. the internal plus the external links. An example about the use of the persistence is reported in Figure \ref{fig:eesempiosemplice}, where two subsets $V_1$ and $V_2$ with different size are considered. Communities are defined as the node subsets with maximum $\alpha(V_C)$, however, one should be careful. By definition, $\alpha(V_C)$ is value in $[0,1]$. The extreme cases refer to the situation in which $V_C$ is a singleton, e.g. $\alpha(V_C)=0$, and $V_C =V$, e.g. $\alpha(V_C)=1$. As a result, the plain optimization of $\alpha(V_C)$ is misleading, as no subset can be better than the whole set $V$. Rather, best values of $\alpha(V_C)$ can be calculated by constraining the cardinality $|V_C| = k$ and then determining the community from the trade-off between $k$ and $\alpha(V_C)$.  

\section{Integer programming problem formulation}
\label{problemformulation}    

Given the graph $G=(V,E)$ previously defined,
the aim of this section is to formulate the maximization of $\alpha(V_C)$, $V_C \subset V, G_C \mbox{ connected}$, as mathematical programming. Problem variables are the boolean $x_i, i = 1,\ldots,n$, taking value 1 if $i \in V_C$, $0$ otherwise. Given the trade-off between size $k$ and $\alpha(V_C)$, the following constraint should be imposed:
    \begin{equation}
        \sum_{i\in V} x_i =  k. 
    \end{equation}
Next, for any pair of nodes $i,j \in V$, let:
    \begin{equation}
    z_{ij} = x_i x_j, \hspace{3mm} \forall i,j \in V,
    \label{zetaij}
    \end{equation}
    and
   \begin{equation}
    w_{ij}= (1-x_i)(1-x_j), \hspace{3mm} \forall i,j \in V,
    \label{doppiavij}
	\end{equation}	
    so that $z_{ij}=1$ if both nodes $i$ and $j$ are in the subset $V_C$, whereas $w_{ij}=1$ if both nodes $i$ and $j$ do not belong to the subset $V_C$. Using variables $w$ and $z$, we have that the persistence probability in formula \ref{formula_alfa} can be written as:
        
    \begin{equation}    
    \alpha(V_C) = \frac{\sum_{i = 1}^{n-1}\sum_{j = i + 1}^n a_{ij}z_{ij}}{\sum_{i = 1}^{n-1}\sum_{j = i + 1}^n a_{ij}(1 - w_{ij})}
    \label{obj:fun}
    \end{equation}

As the persistence must be calculated for connected $V_C$ only, we must impose the connectivity of the subgraph induced by $V_C$. In \cite{Benati2022}, a similar problem is addressed, that is finding the optimal graph-connected clique. Here, we adopt one of the methods proposed there to impose node connectivity: the flow-based approach and, for clarity sake, we briefly describe such method.
    
The subgraph $G_C = (V_C, E_C)$ is connected if a node in $V_C$, denoted as source node, can send a unit of flow to any other node of $V_C$ through an auxiliary digraph $G_{D} = (V,A)$, being $A$ the set of arcs defined in a way that both arcs $(i,j), (j,i) \in A$ if $(i,j) \in E$. 
Variables $f_{ij}$ for all pairs $i,j = 1, \ldots, n$ must be introduced, they will correspond to the flow from node $i$ to node $j$. For a given $V_C$, the source is identified as the node $j$ with maximum index, that is $j = \max\{i|x_i = 1\}$. The flow leaves $j$ to satisfy a demand of one unit flow from all the other nodes of the set $V_C$. To be $G_C$ connected, $V_C$ must allow a feasible solution to these set of constraints: 
     
     \begin{equation}
     	\sum_{j \in V: (i,j) \in A} f_{ij} \le \sum_{j \in V: j <i} z_{ji}+\sum_{j \in V: i<j} z_{ij}, \quad \forall i \in V,
     	\label{upper_bound_flow}
     \end{equation}
	 \begin{equation}
	\sum_{i \in V: (i,p) \in A} f_{ip}- \sum_{i\in V: (p,i) \in A} f_{pi}  \ge x_{p}+(n-2)( x_j-1),
	\,  \forall p,j \in V:  j>p,
	\label{flow_conservation_law}
	 \end{equation}
 	\begin{equation}
 	f_{ij}\ge 0,  \quad \forall (i,j)\in A.
 	\label{fij}
 	\end{equation}
 

Constraint \eqref{upper_bound_flow} states an upper bound on the flow level leaving any node $i \in V$. The right hand side can be positive or null. If a node $i$ does not belong to $V_C$, then $z_{ij}$ and $z_{ji}$ are both zero, therefore no flow can leave the node $i$. On the other hand, if the node $i$ belongs to $V_C$ the right hand size is equal to $|V_C|-1$. Constraint \eqref{flow_conservation_law} expresses a flow conservation law. It is defined for all pairs of nodes $j,p \in V$. 
The left hand side represents the flow net balance for the node $p$, i.e. the difference between entering and the leaving flows. If two nodes $j$ and $p$ belong to the subset $V_C$, the left hand side of constraint \eqref{flow_conservation_law} is at least one. When the subgraph induced by $V_C$ is not connected, the flow leaving $p$ is bounded by $|V_C|-1$ (by \eqref{upper_bound_flow}), but the flow entering in $p$ is less than that leaving $p$, because of the flow leaving (at least) a node in $V_C$ that passes through a node outside $V_C$.
On the contrary, if one, or both, nodes ($j$ and $p$) are not included in subset $V_C$, the right hand side of constraint \eqref{flow_conservation_law} is zero or negative. In these last cases, constraint \eqref{flow_conservation_law} is bounded, by constraint \eqref{upper_bound_flow}, to  $|V_C|-1$. 
Finally, condition \eqref{fij} states a non-negative constraint for the unit of flow.   

Hence, the maximum persistence problem is defined as follows. It is an instance of integer fractional programming, due to equation \eqref{p0:max}, and it is denoted as Problem \textbf{P}.     
    
    \begin{align}
    {\bf P:} 
    \max_{x_{i}} & \displaystyle \hspace{2mm} \frac{\sum_{i = 1}^{n-1}\sum_{j = i + 1}^n a_{ij}z_{ij}}{\sum_{i = 1}^{n-1}\sum_{j = i + 1}^n a_{ij}(1 - w_{ij})} \label{p0:max}\\
     \text{s.t.}  \nonumber\\
    & \label{p0:c1}\hspace{2mm} \sum_{i\in V} x_i =  k, \\
    & \hspace{2mm} z_{ij} = x_i x_j \quad \forall i,j \in V, \label{p0:c2}\\
     & \hspace{2mm} w_{ij} = (1-x_i)(1-x_j) \quad \forall i,j \in V, \label{p0:c3}\\
     & \hspace{2mm} \sum_{j \in V: (i,j) \in A} f_{ij} \le \sum_{j \in V: j <i} z_{ji}+\sum_{j \in V: i<j} z_{ij}, \quad \forall i \in V,\label{p0:f1}\\
     & \hspace{2mm} \sum_{i \in V: (i,p) \in A} f_{ip}- \sum_{i\in V: (p,i) \in A} f_{pi}  \ge x_{p}+(n-2)( x_j-1),
    \,  \forall p,j \in V: j>p, \label{p0:f2}\\
     & \hspace{2mm} f_{ij}\ge 0,  \quad \forall (i,j)\in A\label{p0:f3}\\
    &\hspace{2mm} x_{i}\in \{0,1\}, \quad \forall i\in V.\label{p0:c4}
    \end{align}
    
The objective function in \eqref{p0:max} expresses the persistence probability in terms of variables $z_{ij}$ and $w_{ij}$ defined in \eqref{zetaij} and \eqref{doppiavij}, which, in turn, are formally represented by constraints \eqref{p0:c2} and \eqref{p0:c3}, respectively. \eqref{p0:c1} is the constraint on the cardinality of $V_C$. Constraints \eqref{p0:f1}, \eqref{p0:f2}, and \eqref{p0:f3} ensure that the subgraph $G_C$ is connected. The variable $x_i$ is defined by constraint \eqref{p0:c4}.

The problem $\mathbf{P}$ can be converted into a mixed-integer linear problem  through the Clique Partitioning problem inequalities and the Charnes-Cooper linearization, \citep{charnes1962}.

An auxiliary variable $u$ is introduced such that: 
    
\begin{equation}\label{dnm}
u = \frac{1}{{\sum_{i = 1}^{n-1}\sum_{j = i + 1}^n a_{ij}(1 - w_{ij})}}
\end{equation}

and set $h_{ij} := u\, z_{ij}$ and $l_{ij} := u\, w_{ij}$. 
The Charnes-Cooper linearization allows us to write the objective function in \eqref{p0:max} as
	
	

	

	
	\begin{equation} \label{uzij}
	\sum_{i=1}^{n-1}\sum_{j=i+1}^{n} a_{ij}h_{ij}.
	\end{equation}
	
    Using the Clique Partitioning problem inequalities, variables $z_{ij}$ and $w_{ij}$, are equivalent to:

	$$
	z_{ij} = x_i x_j \Longleftrightarrow 
	\begin{cases}
		z_{ij} \le x_i \\
		z_{ij} \le x_j \\
		z_{ij} \ge x_i + x_j - 1
		\end{cases}
		\hspace{3mm} \forall i,j\in V
	$$
	
		$$
	w_{ij} = (1-x_i)(1-x_j) \Longleftrightarrow
	\begin{cases}
		w_{ij} \le 1-x_i \\
		w_{ij} \le 1-x_j \\
		w_{ij} \ge (1-x_i)+(1-x_j) - 1
	\end{cases}
		\hspace{3mm} \forall i,j\in V.
	$$
	


We fix an appropriate large number $M$ representing the upper bound shared by  $z_{ij}$ and $w_{ij}$. The quadratic terms are linearized as follows

	$$
	h_{ij} = uz_{ij} \Longleftrightarrow 
	\begin{cases}
		h_{ij} \le u \\
		h_{ij} \le Mz_{ij}\\
		h_{ij} \ge u - M(1-z_{ij})
	\end{cases}
		\hspace{3mm} \forall i,j\in V
	$$
	and
		$$
	l_{ij} = uw_{ij} \Longleftrightarrow 
	\begin{cases}
		l_{ij} \le u \\
		l_{ij} \le Mw_{ij}\\
		l_{ij} \ge u - M(1-w_{ij})
	\end{cases}
		\hspace{3mm} \forall i,j\in V.
	$$
	
	Note that if $z_{ij}=1$, the inequalities are satisfied only for $h_{ij}=u$, whereas if $z_{ij}=0$ then the inequalities are satisfied only for $h_{ij}=0$, that is exactly the meaning of $h_{ij} = uz_{ij}$ (similar considerations are valid for $l_{ij}$). 
	
	Hence, the linearized problem can be written as following:
    \begin{align}
    \label{p1:fo}
    {\bf P1:} 
    \max_u & \displaystyle 
    \sum_{i=1}^{n-1}\sum_{j=i+1}^{n} a_{ij}h_{ij}\\
    \text{s.t.} & \nonumber\\
    \label{f1_sum}
   & \sum_{i\in V} x_i =  k.\\
    \label{p1:c1}
    &z_{ij} \le x_i \quad \forall i,j \in V: i<j, \\
    \label{p1:c2}
    &z_{ij} \le x_j \quad \forall i,j \in V: i<j \\
    \label{p1:c3}
    &z_{ij} \ge x_i + x_j - 1 \quad \forall i,j \in V: i<j \\
    \label{p1:c4}
    &w_{ij} \le 1 - x_i \quad \forall i,j \in V: i<j  \\
    \label{p1:c5}
   &w_{ij} \le 1 - x_j \quad \forall i,j \in V: i<j  \\
    \label{p1:c6}
    &w_{ij} \ge (1 - x_i) + (1 - x_j) - 1 \quad \forall i,j \in V: i<j \\
    \label{p1:c7}
    & h_{ij} \le u, \quad \forall i,j \in V: i\ne j\\
    \label{p1:c8}
    & h_{ij} \le M\,z_{ij}, \quad \forall i,j \in V: i\ne j \\
    \label{p1:c9}
    & h_{ij} \ge u - M(1-z_{ij}), \quad \forall i,j \in V: i\ne j \\
    \label{p1:c10}
    & h_{ij} \ge 0 \\
    \label{p1:c11}
    & l_{ij} \le u,\quad \forall i,j \in V: i\ne j\\
    \label{p1:c12}
    & l_{ij} \le M\,w_{ij}, \quad \forall i,j \in V: i\ne j\\
    \label{p1:c13}
    & l_{ij} \ge u - M(1-w_{ij}), \quad \forall i,j \in V: i\ne j \\
    \label{p1:c14}
    & l_{ij} \ge 0, \\
    \label{p1:c15}
    & \sum_{i=1}^{n-1} \sum_{j=i+1}^{n} a_{ij} \cdot (u - l_{ij}) = 1,\\
    \label{f2:0}
    &\sum_{j \in V: (i,j) \in A} f_{ij} \le \sum_{j \in V: j <i} z_{ji}+\sum_{j \in V: i<j} z_{ij}, \quad \forall i \in V,\\
    \label{f2:1}
    &\sum_{i \in V: (i,p) \in A} f_{ip}- \sum_{i\in V: (p,i) \in A} f_{pi}  \ge x_{p}+(n-2)( x_j-1),
    \,  \forall p,j \in V:  j>p, \\
    \label{f2_dom}
    &f_{ij}\ge 0,  \quad \forall (i,j)\in A\\
   \label{p1:boolean}
    &x_{i}\in \{0,1\}, \quad \forall i\in V
    \end{align}

    The objective function in \eqref{p1:fo} represents the persistence probability in terms of the variable $h_{ij}$ defined through the auxiliary variable $u$ in \eqref{dnm}. Constraints \eqref{p1:c1},\eqref{p1:c2}, and \eqref{p1:c3} reply the constraint \eqref{p0:c1} with the Clique Partitioning inequalities of $z_{ij}$.
    The quadratic term $l_{ij}$ is linearized by the constraints \eqref{p1:c11},\eqref{p1:c12}, \eqref{p1:c13}, and \eqref{p1:c14}. Constraint \eqref{p1:c15}
    writes in terms of equality  the auxiliary variable $u$ introduced in \eqref{dnm}.
    Finally, constraints \eqref{p1:c7},\eqref{p1:c8},\eqref{p1:c9}, and \eqref{p1:c10} express the linearization of $h_{ij}$. 

    \bigskip
    
    \section{Heuristic algorithms for the optimal persistent community}
    \label{Heuristic}
    
    In this section we present some heuristic algorithms for finding the community with maximum persistence. The algorithms that we introduce are (in order of computational times): 
    
    \begin{itemize}
        \item A randomized-greedy procedure called Random Shrink (RS).
        \item A merge procedure called Random Shrink Interchange (RSI).
        \item Two long-term heuristic search called Random Shrink Variable Neighborhood Search (RSVNS) and Constrained Random Restart (CRR). 
    \end{itemize}
    
    One practical difficulty of applying the persistent index to community detection is that a user does not know in advance what is the size $k$ of the optimal community. Moreover, as we discussed previously, the persistent index tends to increase as the size of a community increases as well, ranging from 0, when the community is a singleton, to 1, when the community is composed of all nodes of $V$. As we will show in Section \ref{test}, at least in some cases, there is a way to determine $k$, but only after that approximate values of the persistence have been calculated for all $k$  belonging to the set $\{k_l, \dots, k_u\}$. Therefore, the first need of a user is a fast and reliable algorithm to calculate the persistence for communities of size $k$ in a range, and this is the purpose of algorithm Random Shrink.
    Next, after determining $k$, the incumbent solution calculated by Random Shrink can be improved replacing a node of a community with an external one, using the interchange function. When the interchange function cannot improve a solution, we say that the solution is a local optimum, and we call ``basin of attraction" the set of all communities that calculate the same local optimum when interchange is applied. Of course, the interchange subroutine can be applied to many starting solutions to find different local optima, but many strategies are possible and the most effective depends on the problem. We test two of these strategies, one based on variable neighbors and the other on random restart. In variable neighbors search the best solution found is slightly perturbed to escape the local optimum and to continue the interchange in another basin of attraction, that is close to the previous one. If the local optima are close to each other, and the optimal community is similar to (e.g. it overlaps) other communities that are only locally optimal, this strategy can be effective. Conversely, if local optima are distant, so that the optimal solution does not overlap with other local optima, then it is more convenient to forget about them and continue the search though random restart. Hopefully, the new starting solution is a community of a complete different basin of attraction.
    
    \subsection{Random Shrink}
    
    In most application, the researcher does not know in advance the correct size $k$ of the community with the optimal persistence. Rather, as often happens with clustering algorithms, optimization should consider various levels of $k$ before deciding the best one. It may be quite time expensive to engage the algorithms in a thorough calculation of the optimal communities before knowing the exact value of $k$. Instead, it would be convenient to be satisfied with approximate values of $\alpha$ calculated quickly, hence to concentrate the computational resources after that $k$ has been determined.
    
    The Algorithm~\ref{alg:GreedyMerge} that is presented below takes advantage from the fact that, for two non-overlapping
    communities: 1) it takes linear computational time to check whether their union is a connected subset; 2) it takes constant computational time to calculate the resulting value $\alpha$. Hence, starting from a collection of non overlapping communities, they can be progressively merged to find (almost) optimal $\alpha$'s for a whole range of $k$. Moreover, to enhance diversification, the process can be repeated many times with different starting solutions with a cheap computational cost. The method has been used before in \cite{Benati2022}, 
    where it was found to be an effective and reliable tool to obtain a collection of almost optimal communities quickly.
    
    In a general step of the algorithm, nodes of $V$ are partitioned into groups $\cCal = \{C_1,\ldots,C_m\}$, with $m \le n$. Let be  $C_q,C_l \in \cCal$, the following quantities are given:
    \begin{itemize}
    	\item $\mathcal{E}^{in}_q = \#\{\mbox{edges with both end nodes in $C_q$}\}$;
    	\item $ \mathcal{E}^{out}_q = \#\{\mbox{edges with one end node in $C_q$ and the other end node not in $C_q$}\}$;
    	\item $ \mathcal{A}_{ql} = \#\{\mbox{edges with one end node in $C_q$ and the other end node in $C_l$}\}$.
    \end{itemize}
    
    Given those quantities, two clusters $C_q$ and $C_l$ can be merged if $ \mathcal{A}_{ql} \ge 1$, and next it is easy to calculate the coefficient $\alpha(C_{ql})$ of the group $C_{ql} = C_q \cup C_l$:
    
    \begin{equation}
    	\alpha(C_{ql}) = \frac{\mathcal{E}^{in}_q + \mathcal{E}^{in}_l + \mathcal{A}_{ql}}
    	{\mathcal{E}^{in}_q + \mathcal{E}^{in}_l + \mathcal{E}^{out}_q + \mathcal{E}^{out}_l - \mathcal{A}_{ql}}
    \end{equation}
    
So, values $\alpha(C_{ql})$ can be calculated for all $q,l$ pairs and then the best one is selected. Next, the partition $\cCal$ is updated by deleting $C_q, C_l$ from it, but inserting $C_{ql}$ and the process repeated until $|\cCal| = 1$, that is, all subsets are merged. If a data structure containing	$\mathcal{E}^{in}_q, \mathcal{E}^{out}_q, \mathcal{A}_{ql}$ is available from the beginning, then data can be updated in linear time whenever two subsets are merged.
    
The pseudo code of the Algorithm Random Shrink is presented in Algorithm \ref{alg:GreedyMerge}. Merging begins with clusters containing one node, see Line \ref{alg:GreedyMerge:start}. 
While instructions lead the choice of $C_q$ and $C_l$: in first iterations, e.g, $t \le $ {\it max\_random\_step} in Line \ref{alg:GreedyMerge:randomTest}, subsets are chosen at random (Line \ref{alg:GreedyMerge:randomChoice}), in order to diversify the search between different starts to explore different basins of attraction. After that, the algorithm continues in a greedy way, merging clusters that obtain the best local solution $\alpha_k$ (updated in Line \ref{alg:GreedyMerge:update}). When some $\alpha_k$ is 
an optimum, values $\baralpha_k$ are updated (Line \ref{alg:GreedyMerge:max}). Of course, when some $\baralpha_k(C_k)$ is updated, the  corresponding optimal set $C_k$ is updated as well (not reported in the pseudocode). Finally, the process is repeated {\it max\_start} times, see Line \ref{alg:GreedyMerge:restart}, and the fact that the first merging are random could guarantee a sufficient diversification of the search.
    
    \begin{algorithm}
        \caption{Random Shrink}
    	\label{alg:GreedyMerge}
    	\DontPrintSemicolon
    	\KwIn{$G=(V,E)$.}
    	\KwResult{$V_C = \{C^*_2, \ldots, C^*_{n-1}\}$
    		where $C^*_k$ ($2\leq k\leq n -1$) is a $k$-connected subset of $V$.}
    	\SetKwInOut{Parameters}{Parameters}
        \Parameters{$max\_start$ repetition numbers, $max\_random\_step$, number of random choices.}
    	\BlankLine
    	
    	$\baralpha^k = 0$ for $k = 2,\ldots,n-1$\;    
    	\For{$s = 1,\ldots,$max\_start}{\label{alg:GreedyMerge:restart}
    		$\alpha^k = 0$ for $k = 2,\ldots,n-1$\;
    		$\cCal^{0} = \{\{1\},\ldots,\{n\}\}$\;\label{alg:GreedyMerge:start}
    		it = 0 \;
    		\While {$|\cCal^{it}| > 1$}{
    			t = it \;
    			\If{ t $\le$ max\_random\_step}{ \label{alg:GreedyMerge:randomTest}
    				Select randomly $C_q, C_l \in \cCal^t$, such that $G[C_q \cup C_l]$ is connected\;\label{alg:GreedyMerge:randomChoice}
    				$\alpha(C_{ql}) =\alpha(C_q \cup C_l)$\;
    			} \Else {		
    				$\alpha(C_{ql}) = \max\{\alpha(C_q \cup C_l)| C_q, C_l \in \cCal^t$, $G[C_q \cup C_l]$  connected \}\;\label{alg:GreedyMerge:update}
    			} 
    			$\cCal^{t + 1} = (\cCal^{t} - C_q - C_l) \cup C_{ql}\}$ \; 
    			$k = |C_{ql}|$: $\alpha_k = \max\{\alpha(C_{ql}),\alpha_q\}$ \;
    			it = it + 1\;
    			
    		} 
    		$\baralpha^k = \max\{\baralpha^k, \alpha^k\}$, for $k = 2,\ldots,n-1$.\label{alg:GreedyMerge:max}
    	} 
    	
    \end{algorithm}

When the Algorithm Random Shrink concludes, the researcher has at disposal persistence values $\baralpha_k$ for a range of $k$'s, from which to select the best size $k$. How to select $k$ is explained in the computational tests section (Section \ref{test}).

\subsection{Random Shrink Interchange}

After selecting the community size $k$, computational resources can be invested to improve the objective function $\alpha_k$. The Interchange function (Function~\ref{alg:Interchange}) attempts to maximize the $\alpha$-value of a $k$-connected subset $V_C$ of $V$  by replacing one node at a time while keeping the connectivity. The function begins with the initial subset $V_C$ of $V$ containing the candidate nodes and continues to exchange an inner node $h$ ($h\in V_C$) with an outer node $h'$ ($h'\in (V - V_C)$) to obtain a better  $\alpha$-solution   (updated in Line \ref{alg:Interchange:update}). When no more $\alpha$-improvement can be found the function ends, see Line \ref{alg:Interchange:end}, returning the current $k$-connected subset $V_C$.
The computational time of the function can be high, due to the fact that, while it is fast to check whether a candidate entering node is connected to other nodes in $V_C$, a candidate exiting node can leave $V_C$ only if it will not break the connectivity, that must be checked by an appropriate subroutine.

 \begin{function}
    \caption{Interchange()}
    \label{alg:Interchange}
    \DontPrintSemicolon
    \SetKwProg{Fn}{}{}{}
    \SetKwFunction{Interchange}{Interchange}    
    \Fn(){\Interchange{$G$, $V_C$}}{
            \KwIn{$G=(V,E)$, $V_C$ is a $k$-connected subset of $V$.}
            \KwResult{a $k$-connected subset of $V$.}
            \BlankLine
            \While{ \texttt{True}}{
                $\alpha(C^h_{h'}) = \max\{\alpha(C^j_i)| \forall i\in V_C, j\in(V - V_C) \text{ and } C^j_i  \text{ connected}\}$\;\label{alg:Interchange:update}
                \eIf{$\alpha(V_C) < \alpha(C^h_{h'})$}{
                    remove $h$ from $V_C$ and add $h'$ to $V_C$\;  
                }
                {
                    \Return $V_C$\;\label{alg:Interchange:end}
                }
            }
        }
\end{function}

 We combine the algorithm Random Shrink with the Interchange function obtaining a new algorithm that we will call Random Shrink Interchange described in Algorithm~\ref{alg:GreedyMergeInterchange}.
 
 \begin{algorithm}
        \caption{Random Shrink Interchange}
        \label{alg:GreedyMergeInterchange}    
        \SetAlgoLined
        \DontPrintSemicolon
    	\KwIn{$G=(V,E)$, $k$  community size.}
    	\KwResult{$C^*_k$ ($2\leq k\leq n -1$) is a $k$-connected subset of $V$.}
    	\SetKwInOut{Parameters}{Parameters}
        \Parameters{$max\_start$ repetition numbers, $max\_random\_step$, number of random choices.}
    	\BlankLine
    	Let $C_k$ be the $k$-connected subset of $V$ resulting from the Algorithm~\ref{alg:GreedyMerge} applied to $G$\;
    	$C^*_k = \Interchange(G, C_k)$\;\label{alg:GreedyMergeInterchange:optimized}
    \end{algorithm}
 
\subsection{Random Shrink Variable Neighborhood Search}

After the application of the Interchange function, a first 
optimal solution $V_C \subset V$, $|V_C| = k$ has been determined.  
The Tree Variable Neighborhood Search function (Function~\ref{alg:vns}) attempts to improve the $\alpha$-value of $V_C$  by replacing a subset $\rCal$ of $V_C$ with a random subset of nodes $\iCal$ taken from the neighbors of $V_C$ and obtaining a new solution $V_C' = (V_C - \rCal) \cup \iCal$. Next, the Interchange function is applied to $V_C'$, hopefully to find a new local optimum. If $V_C$ is almost optimal, this way of diversifying the search can be effective. 

Calculating $V_C'$ from $V_C$ can be computationally demanding as $V_C' = (V_C - \rCal) \cup \iCal$ must be connected, that is the reason why we implement a subroutine that exploits the connectivity of a random spanning tree of $G_C$. Specifically, $\rCal$ is obtained in the following way. First, a random node $v\in V_C$ (Line \ref{alg:vns:random}) is selected and a random spanning tree $T=(V_C,E_T[G_C])$ is constructed from the root $v$. 
Next, we determine the set $\lCal$ of nodes that are leaves of $T$ and select a random number $h$ in the range $h \in \{2,\dots,|\lCal|\}$. Finally, the exiting nodes are a random set $\rCal \subseteq \lCal$, such that $|\rCal| = h$. The way in $\rCal$ is selected guarantees that the subgraph having vertices set $V_C - \rCal$ and edges set $E[V_C - \rCal]$ is connected. The subset $\iCal \subset V$ such that $|\iCal| = h$ is selected at random from the neighboring nodes (different from those belonging to the set $\rCal$) of $V_C - \rCal$. 
The subset of $V$ so obtained is optimized with the Interchange function (Line \ref{alg:vns:optimized}).
Finally, the process is repeated $max\_start$ times to obtain different starting solutions $V_C'$, see Line \ref{alg:vns:restart}.
 
 \begin{function}
        \caption{Tree Variable Neighborhood Search()}
        \label{alg:vns}
        \DontPrintSemicolon
        \SetKwProg{Fn}{}{}{}
        \SetKwFunction{VNS}{VNS}    
        \Fn(){\VNS{$G$, $V_C$}}{
            \KwIn{$G=(V,E)$, $V_C$ is a $k$-connected subset of $V$.}
            \KwResult{a $k$-connected subset of $V$.}
            \SetKwInOut{Parameters}{Parameters}
            \Parameters{$max\_start$ repetition numbers.}
            \BlankLine
    	    \For{$s = 1,\ldots,max\_start$}{\label{alg:vns:restart}
    	        select randomly a node $v\in V_C$\;\label{alg:vns:random}
    	        select randomly a subset $\rCal\subset V$ build from the leaf of a random spanning tree of $v$\;\label{alg:vns:removed}
    	        select randomly a subset $\iCal$ from the neighbors of $V_C-\rCal$ of the same size of $\rCal$\;\label{alg:vns:inserted}
    	        $\overline{V_C} = \Interchange(G, (V_C - \rCal)\cup \iCal)$\;\label{alg:vns:optimized}
                \If{$\alpha(V_C) < \alpha(\overline{V_C})$}{
                    replace  $V_C$ with $\overline{V_C}$\;  
                }
            }
            \Return $V_C$\;\label{alg:Interchange:return}
        }
    \end{function}
 
 We combine the algorithm Random Shrink Interchange with the Tree Variable Neighborhood Search function obtaining a new algorithm that we will call Random Shrink Variable Neighborhood Search described in Algorithm~\ref{alg:GreedyMergevns}.

    \begin{algorithm}
        \caption{Random Shrink Variable Neighborhood Search}
        \label{alg:GreedyMergevns}    
        \SetAlgoLined
        \DontPrintSemicolon
    	\KwIn{$G=(V,E)$, $k$  community size.}
    	\KwResult{$C^*_k$ ($2\leq k\leq n -1$) is a $k$-connected subset of $V$.}
    	\SetKwInOut{Parameters}{Parameters}
        \Parameters{$max\_start\_greedy$ repetition numbers greedy alg , $max\_random\_step$, number of random choices, $max\_start\_vns$ repetition numbers vns function .}
    	\BlankLine
    	Let $C_k$ be the $k$-connected subset of $V$ resulting from the Algorithm~\ref{alg:GreedyMergeInterchange} applied to $G$\;
    	$C^*_k = \VNS(G, C_k)$\;\label{alg:GreedyMergevns:optimized}
    \end{algorithm}

    \subsection{Constrained Random Restart}
    Even though Random Shrink Variable Neighborhood Search algorithm is a flexible tool to explore local optima that are close to the incumbent solution, still it could be the case that the global optimum is much farther away and it can be detected only selecting a completely different starting solution. As the Algorithm~\ref{alg:GreedyMergevns},  Constrained Random Restart algorithm tries to improve a set of starting solutions through the Interchange function, but in this new algorithm the starting solutions are selected at random. Still, as it is important to explore different basins of attraction, we have included a distance constraint about how new starting solutions are selected. First, a new starting node $c$ must have a distance from the previously selected nodes (set $\pCal$) of at least $min\_distance$ (Line~\ref{alg:RandomRestart:test}) and then a new community is built at random around $c$. When, due to distance constraints, $c$ cannot be determined, a new $V_C$ is determined from scratch and the process is repeated until a number of $max\_start$ solutions is attempted (Algorithm~\ref{alg:RandomRestart}).
      
    \begin{algorithm}
        \caption{Constrained Random Restart}
        \label{alg:RandomRestart}    
        \SetAlgoLined
        \DontPrintSemicolon
        \KwIn{$G=(V,E)$, $k$  community size.}
        \KwResult{$V_C^b$ a $k$-connected subset of $V$.}
        \SetKwInOut{Parameters}{Parameters}
        \Parameters{$\dCal$ distance function, $min\_distance$ minimum distance between the starting node, $max\_start$ repetition numbers.}
        \SetKwRepeat{Do}{do}{while}
        \BlankLine
        $it = 0$ \;
    	\While {$it < max\_start$}{
            randomly select a node $c\in V$\;\label{alg:RandomRestart:start:begin}
            $\pCal=\{c\}$\;
            randomly built a $k$-connected subset $V_C$ of $V$\;
            $V_C^b = \Interchange(G, V_C)$\;\label{alg:RandomRestart:start:end}
            \Do{$it < max\_start$}
            {\label{alg:RandomRestart:iteration}
                randomly select a node $c\in V$ such that $\dCal(\pCal, c) \geq min\_distance$\;
                \eIf{$c$ is found\label{alg:RandomRestart:test}} 
                {
                    $it = it + 1$\;
                    add $c$ to $\pCal$\;
                    randomly built a $k$-connected subset $V_C$ of $V$\;
                    $\overline{V_C} = \Interchange(G, V_C)$\;
                    \If{$\alpha(\overline{V_C}) > \alpha(V_C^b)$} {
                        $V_C^b = \overline{V_C}$
                    }
                }
                {
                    \textbf{break}\;
                }
            }
        }
    \end{algorithm}
      
\section{Computational experiments}\label{test}
   
In this section, we test the proposed algorithms of Section \ref{Heuristic} on a family of simulated networks that have been generated according to the methodology proposed in \cite{lancichinetti2008}. That procedure generates synthetic networks that are as close as possible to real networks, which are often characterized by a 
high variability in the nodes degree. We simulate networks of $n$ nodes 
so that each node degree is a random value taken from 
a power law distribution with parameter $\gamma$, minimum and maximum degree $k_{min}$ and $k_{max}$, respectively, and average degree $\langle k \rangle$.
Then, nodes of the graph are partitioned into communities, using the mixing parameter $\mu \in (0,1]$. This parameter represents the fraction of edges that starting from a given node points to nodes outside the community. Conversely, the complement $1-\mu$ is the fraction of edges outgoing from a node and pointing to nodes inside the community. The size of each community is a random value taken from a power law distribution with parameter $\beta$, and it ranges between the minimum $s_{min}$ and the maximum $s_{max}$. The procedure partitions nodes of the network with each node being assigned to only one community.
We simulate $N = 1000$ graphs of size $n \in \{20, 25, 30, 40, 50, 100, 150, 200\}$. For each graph, we set the mixing parameter $\mu = 10\%$, the average node degree $\langle k \rangle = 30\%$, and $s_{min}$ and $s_{max}$ equal to $20\%$ and $50\%$, respectively, parameters $\gamma$ and $\beta$ are set to $2$ and $1$, respectively, that are the lowest values of the intervals indicated in \cite{lancichinetti2008}. 
    
All algorithms have been implemented in \texttt{C++} language. The exact solution is calculated using the solver GuRoBi described in~\cite{gurobi2022}. The simulations have been performed on an iMac Pro 3.2 GHz 8-Core Intel Xeon W with 32 GB of ram. 
    
Our first preliminary test consists in applying the model {\bf P1}, solved with GuRoBi, to simulated networks of size $n$ ranging from 20 to 40. For these particular class of graphs, we set $k = n/2$. In Table \ref{exact}, we reported the average of the computational times on 20 instances (i.e. 20 simulated graphs). As it appears evident, only instances of small size can be solved by Integer Programming: for $n= 30$ times are less than 10 minutes, but they become more than 10 hours for $n=40.$
    
\subsection{Determining the correct community size \texttt{k}}
    
\begin{figure}[htp]
        \centering
        \includegraphics[width=15cm,height=8cm]{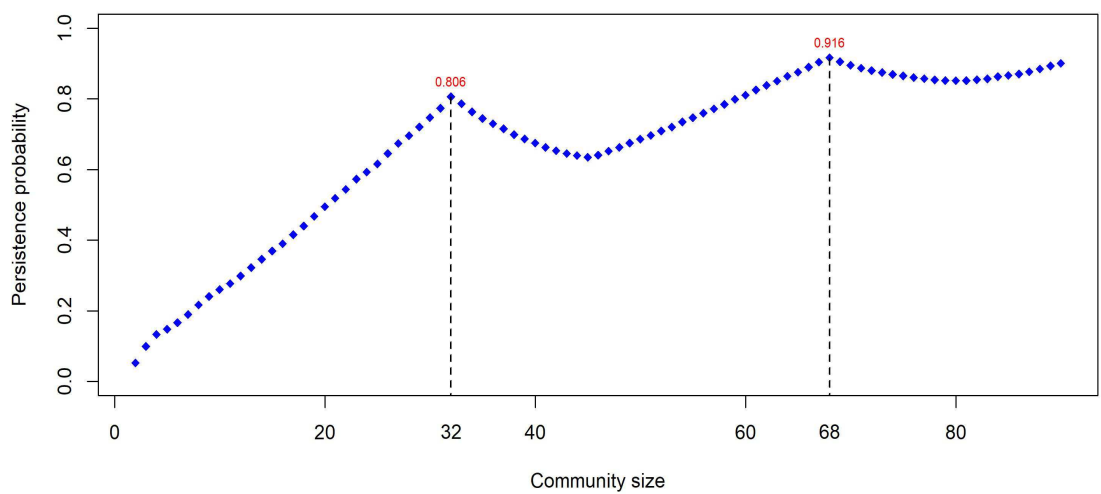}
        \caption{Example of persistence curve.}
        \label{fig:picchi}
    \end{figure}

As we discussed earlier, the computation of the optimal persistence is flawed by the mere fact that the index tends to increase just as the cardinality of the communities increases. In particular, the value of $\alpha$ is close to zero, when communities are small and connected with many other nodes outside the community, and it approaches to one, when communities are composed of almost all the nodes of the network connected mostly with nodes inside the community. An example is reported in Figure \ref{fig:picchi}, that depicts, for a simulated graph, the curve of the persistence $\alpha_k$ as the community size $k$ increases. As it can be seen, the 
trend of the curve is 
almost increasing,
providing an evidence that the global maximum of the function cannot be the unique criterion to select communities. Nevertheless, a closer inspection of the figure reveals that there are some local maxima determined for intermediate value of $k$:  as better motivated later, we guess that they are the correct values of $k$ that determine communities. 
    Therefore, it is important relying on a fast and accurate method for drawing histograms as in Figure \ref{fig:picchi}, that we will call \textit{persistence curve}.
    
    \begin{table}[htp]
    \centering
\begin{tabular}{|l|c|l|l|}
\hline
{\color[HTML]{495057} \textit{n}} & {\color[HTML]{495057} \textbf{num. networks}} & {\color[HTML]{495057} \textbf{$k$}} & {\color[HTML]{495057} \textbf{time/sec}} \\ 
\hline \hline
{\color[HTML]{495057} 20}         & {\color[HTML]{495057} 20}                     & {\color[HTML]{495057} 10}                      & {\color[HTML]{495057} 16.39}             \\ 
{\color[HTML]{495057} 25}         & {\color[HTML]{495057} 20}                     & {\color[HTML]{495057} 13}                      & {\color[HTML]{495057} 48.68}             \\ 
{\color[HTML]{495057} 30}         & {\color[HTML]{495057} 20}                     & {\color[HTML]{495057} 15}                      & {\color[HTML]{495057} 456.10}            \\ 
{\color[HTML]{495057} 40}         & {\color[HTML]{495057} 20}                     & {\color[HTML]{495057} 20}                      & {\color[HTML]{495057} 21040.69}          \\ \hline
\end{tabular}
\caption{Time computation for solving the exact problem }
\label{exact}
\end{table}

Algorithm \ref{alg:GreedyMerge} has been developed for the purpose. It has been tested for graphs of size ranging from $n = 20$ to $n=200$. As already pointed out, for each size, 1000 graphs have been generated and the algorithm is run, determining an approximate value of $\alpha_k$ for all $k \le n$.
We set the number \textit{maxit} equal to 100, 1000 and 10000 starting solutions and we summarize computational results in Table \ref{tbl_greedy}. The left side of the table reports the objective function (i.e., the persistence). It is measured as follows: for each instance (and depending on the size $n$), we calculate the index\footnote{In the sum we exclude the extreme cases $k=1$ and $k=n$ that correspond to trivial cases.} $f_{n} = \sum_{k = 2}^{n-1} \alpha_k$, then, we computed the average of $f_n$ on all the instances for each value of $n$. 

The two central columns report the probability of improvement when the number of starting solutions increases from 100 to 1000 and from 100 to 10000. It is measured, for each instance of size $n$, as the relative frequency in which $\alpha_k$ computed for $\textit{maxit} = 1000$ (or for $\textit{maxit} = 10000$, respectively) is greater than $\alpha_k$ computed for $\textit{maxit} = 100$. Then, we calculate the average on all the instances for each value of $n$.
It can be noticed that, when \textit{maxit} increases, the corresponding increase of the persistence is quite marginal, always less than 1$\%$. However, these results have to be read in the light of the probability of improvement. This probability significantly increases as $n$ grows, for instance passing from  $4 \%$ for $n=20$ to $58 \%$ for $n=200$ in the case in which \textit{maxit} passes from 100 to 1000. 
 
Computational times are reported as average values in the right side of the table. Even though they linearly increase, they could be severe when more starting solution are allowed. For example, considering $n = 200$, solutions with $maxit = 100$ has been obtained in 2 minutes, with $maxit = 1000$ in 20 minutes, with $maxit = 10000$ in more than 3 hours (note that times are multiplied by 10). However, the purpose of Algorithm \ref{alg:GreedyMerge} is to obtain a quick histogram of the $\alpha_k$ values, and the table reveals that reliable data can be obtained with $maxit = 100$.
    
    \begin{table}[htp]
    \centering
    \begin{tabular}{|l|c|c|c|c|c|c|c|c|}
    \hline
    \multicolumn{1}{|c|}{} 
    & \multicolumn{3}{c|}{Persistence} & \multicolumn{2}{c|}{Probab.} & \multicolumn{3}{c|}{Times} \\
    \hline
    \diaghead{\theadfont normalsize}{$n$}{$maxit$} & 100 & 1000 & 10000 & 1000 & 10000 & 100 & 1000 & 10000 \\ 
        \hline \hline
        20 & 11.79 & 11.82 & 11.82 & 0.04 & 0.04 & 0.05 & 0.64 & 6.61 \\ 
        25 & 14.95 & 15.01 & 15.02 & 0.07 & 0.07 & 0.13 & 1.32 & 14.81 \\ 
        30 & 18.42 & 18.51 & 18.52 & 0.10 & 0.11 & 0.26 & 2.68 & 29.52 \\ 
        40 & 23.97 & 24.11 & 24.13 & 0.15 & 0.17 & 0.83 & 8.46 & 87.18 \\ 
        50 & 29.14 & 29.32 & 29.35 & 0.20 & 0.24 & 2.01 & 20.30 & 207.10 \\  
        100 & 54.94 & 55.27 & 55.50 & 0.38 & 0.47 & 31.25 & 315.24 & 3175.45 \\ 
        150 & 81.52 & 82.10 & 82.40 & 0.49 & 0.60 & 152.96 & 1540.27 & 15277.62 \\ 
        200 & 105.93 & 106.80 & 107.24 & 0.58 & 0.69 & 480.37 & 4880.58 & 47114.58  \\
        \hline
    \end{tabular}
     \caption{Persistence mean and computational times Algorithm \ref{alg:GreedyMerge} depending on graph size.}
     \label{tbl_greedy}
    \end{table}
    
    \begin{table}[htp]
    \centering
    \begin{tabular}{|lllll|}
    \hline
        $n$ & p.k first & p.k median & p.k atleast & p.k all \\ \hline \hline
        20 & 0.812 & 0.735 & 0.948 & 0.465\\ 
        25 & 0.771 & 0.665 & 0.956 & 0.426 \\ 
        30 & 0.716 & 0.606 & 0.949 & 0.470 \\
        40 & 0.824 & 0.711 & 0.988 & 0.602 \\
        50 & 0.843 & 0.814 & 0.997 & 0.673 \\
        100 & 0.866 & 0.860 & 1 & 0.754 \\
        150 & 0.823 & 0.846 & 1 & 0.757 \\ 
        200 & 0.777 & 0.779 & 1 & 0.682 \\ \hline
    \end{tabular}
    \caption{Probabilities of finding a correct value of community size $k$.}
    \label{t:size}
\end{table}

From the analysis of a persistence curve, as the one reported in Figure \ref{fig:picchi}, a researcher can guess the correct community size $k^*$ as the one corresponding to a local maximum. Actually, the persistence curve can contain more than one peak, because the networks can contain more than one community or by mere numerical reasons.  In practical applications, we guess that some further substantive analyses may be carried on the communities corresponding to each peak to establish whether they are realistic clusters. Nevertheless, some indications can be inferred by using some automatic procedure. In our tests, we apply two selection rules.  Suppose that $\{k_1^*,\ldots,k_l^*\}$ 
is the set of the sizes corresponding to the local maxima of the persistence curve. The first rule is taking the smallest value $k_1^*$, the second one is taking the median value $k_m^*$ (with $ m = \lfloor (l+1)/2 \rfloor$). The experimental graphs simulated by the algorithm of  \cite{lancichinetti2008} usually contain more than one community, so let $\{k_1,\ldots,k_r\}$ be the set of sizes of the $r$ simulated communities. Therefore, to check the effectiveness of our procedure we control whether $k^*_i \in \{k_1,\ldots,k_r\}$ with $i = 1$ or $m$. Computational results are described in Table \ref{t:size}. There, in the first two columns and for the 1000 graphs generated by each network size $n$, we have reported the relative frequency with which it has been observed $k^*_i \in \{k_1,\ldots,k_r\}$. As it can be seen, for networks of small dimension,  $k^*_1$ is better than $k^*_m$, but the relation reverses for the largest networks. Nevertheless, and in both cases, they are correct guesses for the large majority of the instances. In more than 10\% of the cases, $k_i^*$, $i = 1$ or $m$ are a wrong prediction, but 
further substantive analyses show that (third column) almost always at least one  $k_i, i = 1,\ldots,r$ is between the guessed ones $k_i^*, i = 1, \ldots, l$.  
Moreover, in the fourth column we report how many times all values $k_i, i = 1,\ldots,r$ are contained in the guessed set $\{k_1^*,\ldots,k_l^*\}$ and
this probability significantly increases when the size $n$ grows. To summarize, most of the selected values $k_i^*$ corresponds to true community sizes $k_i$ and our algorithm can reveal them.
    
\subsection{Interchange heuristics}
    
Our final tests regard the computational analysis of the interchange heuristics. Next algorithm take as input the subset size $k^*$ and the subset $V_C$ determined by algorithm Random Shrink, then $V_C$ is attempted to be improved first by the Interchange function, next by perturbing the optimal solution to restart the interchange from different initial clusters. Tested procedures require a minimal amount of parameters: in the Tree Variable Neighborhood Search function, the number of nodes that are replaced from the incumbent solution is a random number between 2 and $k$. The numbers of initial solution that are tested by Algorithm~\ref{alg:GreedyMergevns} and Algorithm~\ref{alg:RandomRestart} is 100.  In Table~\ref{tabella_confronti} it can be seen the comparison between the three heuristics. Results are average values on 1000 test problems, times included the use of Random Shrink. In the first group of columns data about the plain Random Shrink Interchange are reported. In the first column (P.BtRS: Probability Better than Random Shrink), we report the relative frequency in which the Random Shrink Interchange could improve the Random Shrink solution; it can be seen that it happens, but less frequently than what expected. In the smallest sized problems, only 4\% of the Random Shrink solutions were improved by the Random Shrink Interchange. Frequency of improved solutions increases with problem size; nevertheless, only 16\% of the times the Random Shrink has been increased when the network size is $n=200$. As whole, these data points that the Random Shrink is an effective heuristic. The second column (M.diff: mean difference), reports the relative increase of the objective function calculated only for the cases in which an improvement actually occurs. It can be seen that the improvement is more substantial on the smallest problem size than for the largest, indicating that there are cases in which the Random Shrink fails to find the optimal solution and that the improvement can be substantial. The last column (Time) reports the computational times, in which it can be seen that they are rather negligible, as the largest instances are solved in a few more than one second.
    
We compare the restricted diversification, e.g., Random Shrink Variable Neighborhood Search, in which starting solutions are generated close to the local optimum, with the free diversification, e.g. Constrained Random Restart, in which starting solutions are generated through consideration about their distance from previous analyzed regions. Looking at the frequency in which the Random Shrink solution has been improved, we can see that data are in favor of Constrained Random Restart, as they are improved from 9\% of the times when the network size is $50$ to 22\% when the network size is 200 with respect to 6\% to 18\%. Solutions quality is better as well, improving of some 13\% in all problem size with respect to less 10\%. The only comparison in favor of Random Shrink Variable Neighborhood Search is about times that are some half the ones of Constrained Random Restart.

This is due to the fact that starting solutions generated closer to a local optimum are close to an other local optimum as well, but this implies that the diversification strategy is less effective. In summary, data are suggesting that the Constrained Random Restart is the most effective method to improve the results of the Random Shrink. 
    
     \begin{table}[htp]
    \centering
    \begin{tabular}{|l|ccc|ccc|ccc|}
    \hline
    $n$ &

  \multicolumn{3}{c|}{RSI \ref{alg:GreedyMergeInterchange}} &
  \multicolumn{3}{c|}{RSVNS \ref{alg:GreedyMergevns}} &
  \multicolumn{3}{c|}{CRR  \ref{alg:RandomRestart}} 
  
         \\ \hline
        & P.BtRS     &   M.diff  &   Time    &   P.BtRS   &   M.diff  &   Time    &   P.BtRS   &   M.diff  & Time \\ \hline \hline
       50 & 	0.04    &	0.07    &	0.02  &	0.06	&   0.10    &   0.05  & 0.09	&   0.12    &	0.31    
\\ \hline
       100 & 	0.10	&	0.03	&	0.15 &	0.13	&	0.04	&	0.72	& 0.16	&	0.14	&	2.85	
 \\ \hline
       150 & 	0.14	&	0.02	&	0.46 &	0.16	&	0.04	&	2.69	& 0.20	&	0.13	&	10.02	
 \\ \hline
       200 & 	0.16	&	0.01	&	1.14 &	0.18	&	0.03	&	11.30	& 0.22	&	0.13	&	28.60	
 \\ \hline
    \end{tabular}
    \caption{Comparison between the number of times that a combined algorithm found a better solution w.r.t. the Random Shrink algorithm.}
    \label{tabella_confronti}
\end{table}

\section{Application to real data}
\label{real_analysis}

\begin{figure}[htp]
    \centering
    \begin{subfigure}[h]{0.4\textwidth}
        \includegraphics[width=\textwidth]{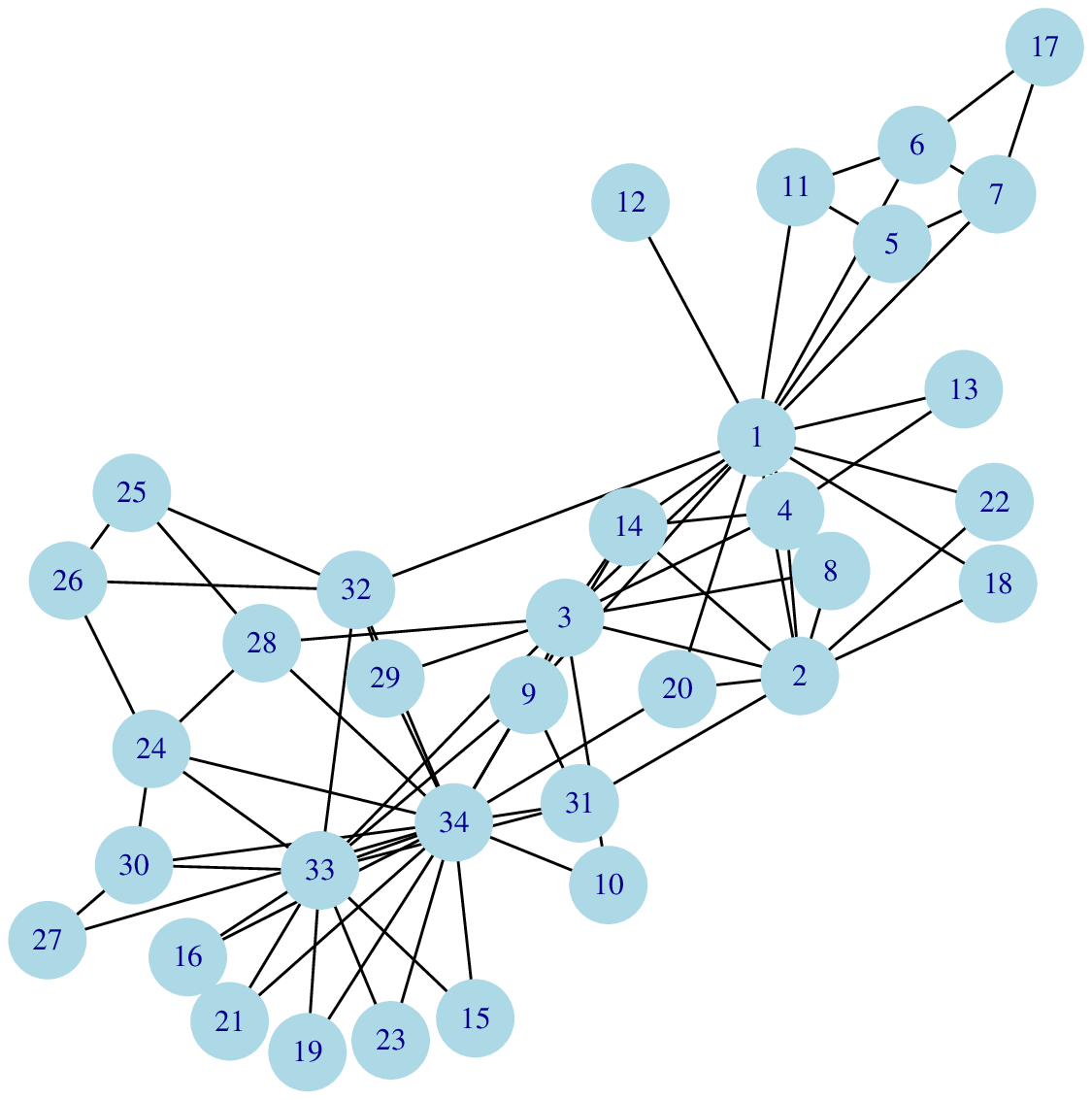}
        \caption{Karate Club network.}
        \label{Karate_club_graph}
    \end{subfigure}
    \begin{subfigure}[h]{0.4\textwidth}
        \includegraphics[width=\textwidth]{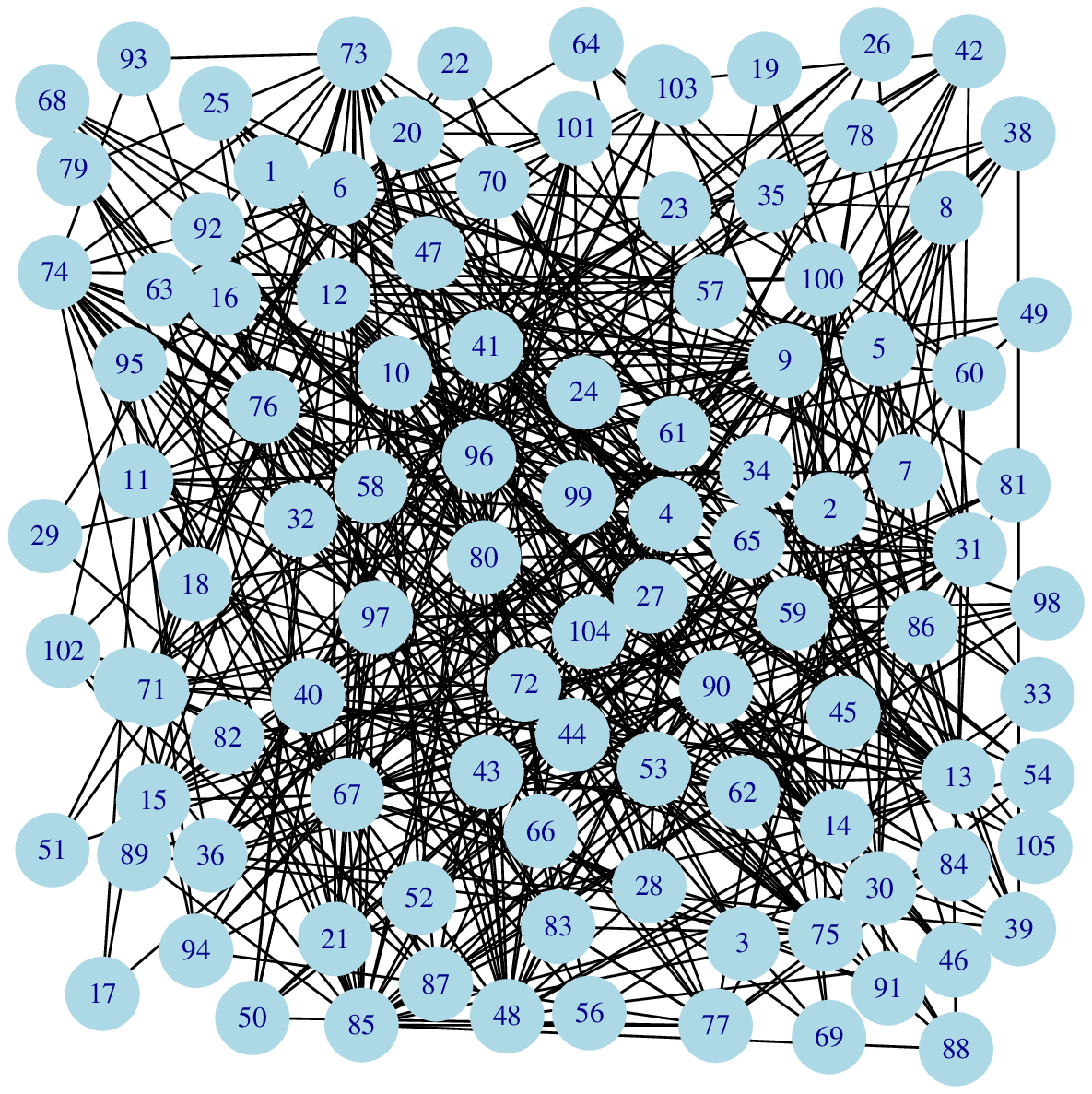}
        \caption{U.S. Political books network.}
        \label{political_books_graph}
    \end{subfigure}
    \caption{Network topologies.}
    \label{Network topologies}
\end{figure}

In this section, we make some computational tests on two simple real networks to verify the reliability of the persistence $\alpha$ to real applications, that is, if it can identify homogeneous groups of nodes interpreted as communities. Considered networks are the Zachary Karate Club (\cite{zachary1977}), and the political books (\cite{krebs2004}).  The Zakary karate club is the network of friendships between the members of a club is a US university, while the political books is the network of co-purchased books about US politics published around year 2004 and sold online by Amazon.com. Figure \ref{Network topologies} represents the topology of these networks. 
    
The two networks are quite different in terms of number of nodes, size, average degree and density. Table \ref{network_features} reports the basic characteristics of the networks. 
    
    \begin{table}[htp]
    \centering
    \begin{tabular}{|lcccc|}
    \hline
        Network & N. of nodes & N. of edges & Average degree & Density \\ \hline \hline
        Karate Club & 34 & 78  & 4.58 & 0.14  \\ 
        Political books & 105 & 441 & 8.4 & 0.081 \\ \hline
    \end{tabular}
    \caption{Network characteristics.}
    \label{network_features}
\end{table}

We apply our methodology, e.g. Algorithm \ref{alg:GreedyMerge}, to compute the peaks of the persistence probability $\alpha$ and then improving those results with Constrained Random Restart, and finally we compare its results with two well-known community detection methods, i.e. Walktrap and Louvain method (\citep{pons2006} and \citep{Blondel2008}). Similar to the idea on which the persistence is based, the Walktrap method assumes that a random walk on a graph tends to remain inside a community. Conversely to the persistence, the community is detected using specific structural distance between vertices and then a hierarchical clustering algorithm is applied. Again, by analogy with the persistence, the Louvain method is based on the maximisation of an index, the modularity score, but the community is detected using spectral decomposition methods and not by optimization. These latter methods are designed to make community detection, that is, finding a partition of nodes, and therefore their output is composed of many communities, while our method is designed to find just one community. Nevertheless, we can discuss the consistency of the results form the three methods.  

\begin{figure}[htp]
	\begin{center}
	    \subfloat[Persistence probability]{\includegraphics[width=1\textwidth]{{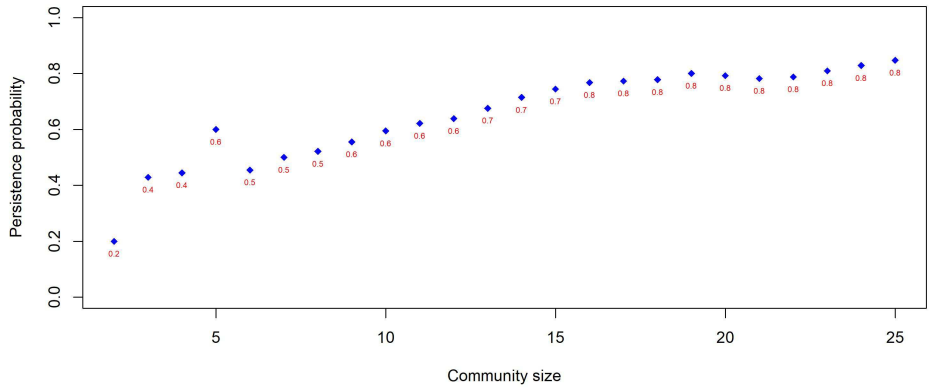}}
	    \label{KC_persistence}
	    }
	\end{center}
	\begin{center}
	    \subfloat[Walktrap algorithm]{\includegraphics[width=0.5\textwidth]{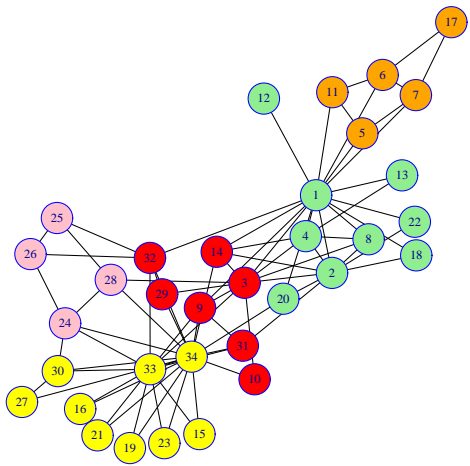}
	    \label{KC_walktrap} \hspace{1pt}%
 	    }
 	    \subfloat[Louvain methodology]{\includegraphics[width=0.5\textwidth]{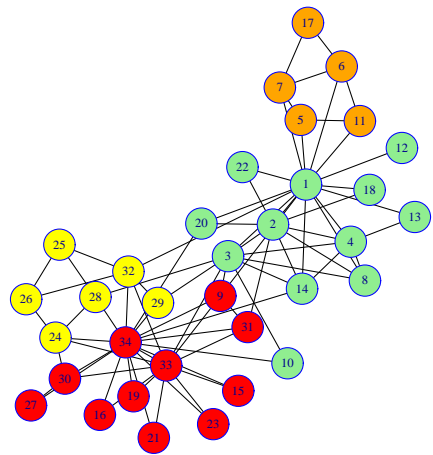}
 	    \label{KC_louvain}
 	    }
	\end{center} 	 
	\caption{Comparison of community detection methods and persistence probability for Zachary karate club network.}
\end{figure}

Figure \ref{KC_persistence} reports the persistence probability ($y$-axis) varying $k$ ($x$-axis) for the Zachary karate club network (figure \ref{Karate_club_graph}). The curve shows that there is an evident peak of the persistence probability 
for $k=5$ corresponding to subset $V_C = \{5,6,7,11,17\},$ in which $\alpha(V_C)=0.60$. This value is well above the persistence of communities of size 4 and 6 showing that those nodes form a well-separated cluster. The  Walktrap and the Louvain method (see Figures \ref{KC_walktrap} and \ref{KC_louvain}), detect this community as well: it is the group with orange color, a specific peripheral group of friends densely connected, but with few links with the rest of graph. Next, the persistence curve reveals a second local maximum, corresponding to the subset $V_C' = \{ 3, 9, 10, 15, 16, 19, 21, 23, 24, 25, 26, 27, 28, 29, 30, 31, 32, 33, 34 \}$ of size $k' = 19$. With respect to the Walktrap communities, $V_C'$ is the union of all the yellow nodes, all the pink nodes, all the red nodes except node 14; with respect to the Louvain communities, $V_C'$ corresponds to all the red nodes, all the yellow nodes, and moreover the green nodes 3 and 10. 
It is worth noting that we also found a group of 10 nodes ($\tilde k=10 = n-(k+k')$ in the persistence curve) that, even though it does not correspond to a local peak, it is composed of the remaining nodes of the green community, that is, excluding nodes 3 and 10. 
In conclusion, the persistence index revealed communities similar to the other methods, possibly allowing some aggregation and with some peripheral nodes resolved differently. Moreover in this case, even though the method purpose is not finding a partition, still the outcome can be interpreted as such.  

\begin{figure}[htp]
	\begin{center}
	    \subfloat[Persistence probability]{\includegraphics[width=1\textwidth]{{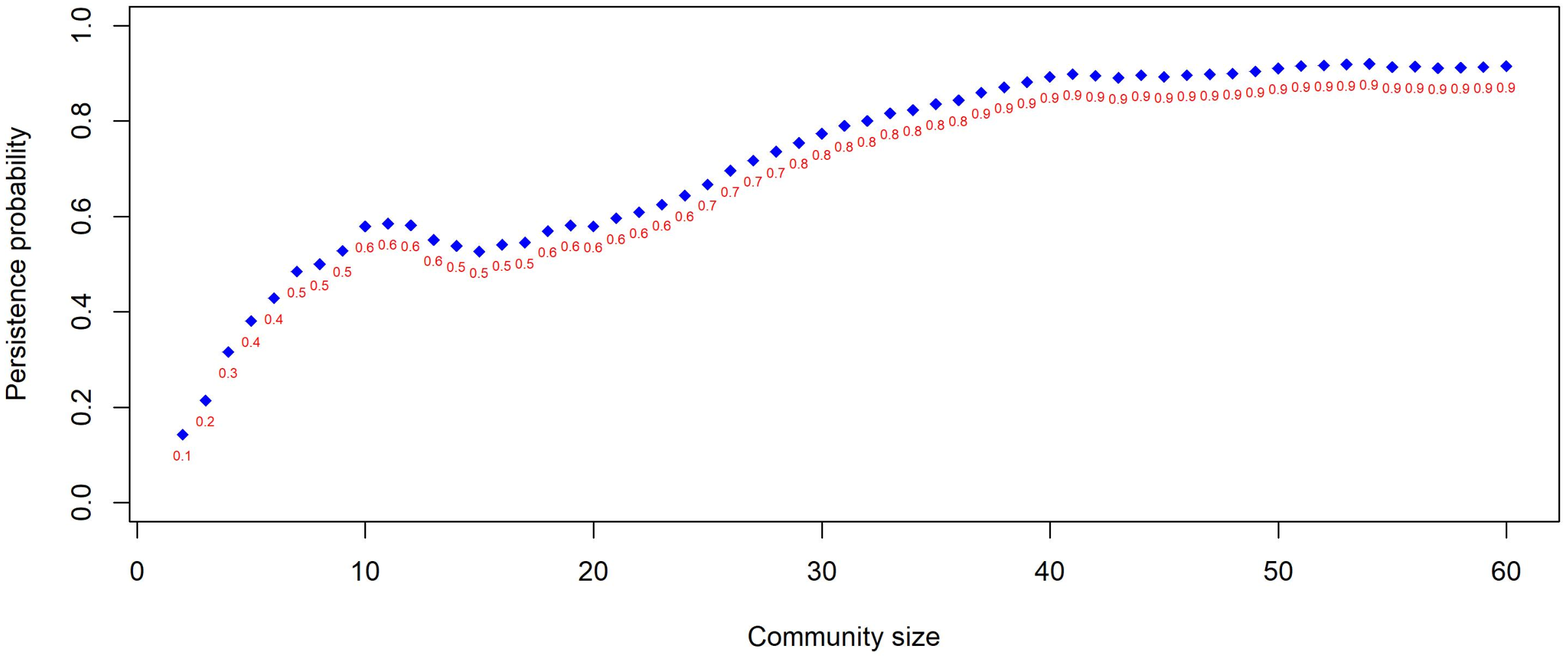}}
	    \label{Book_persistence}
	    }
	\end{center}
	\begin{center}
	    \subfloat[Walktrap algorithm]{\includegraphics[width=0.5\textwidth]{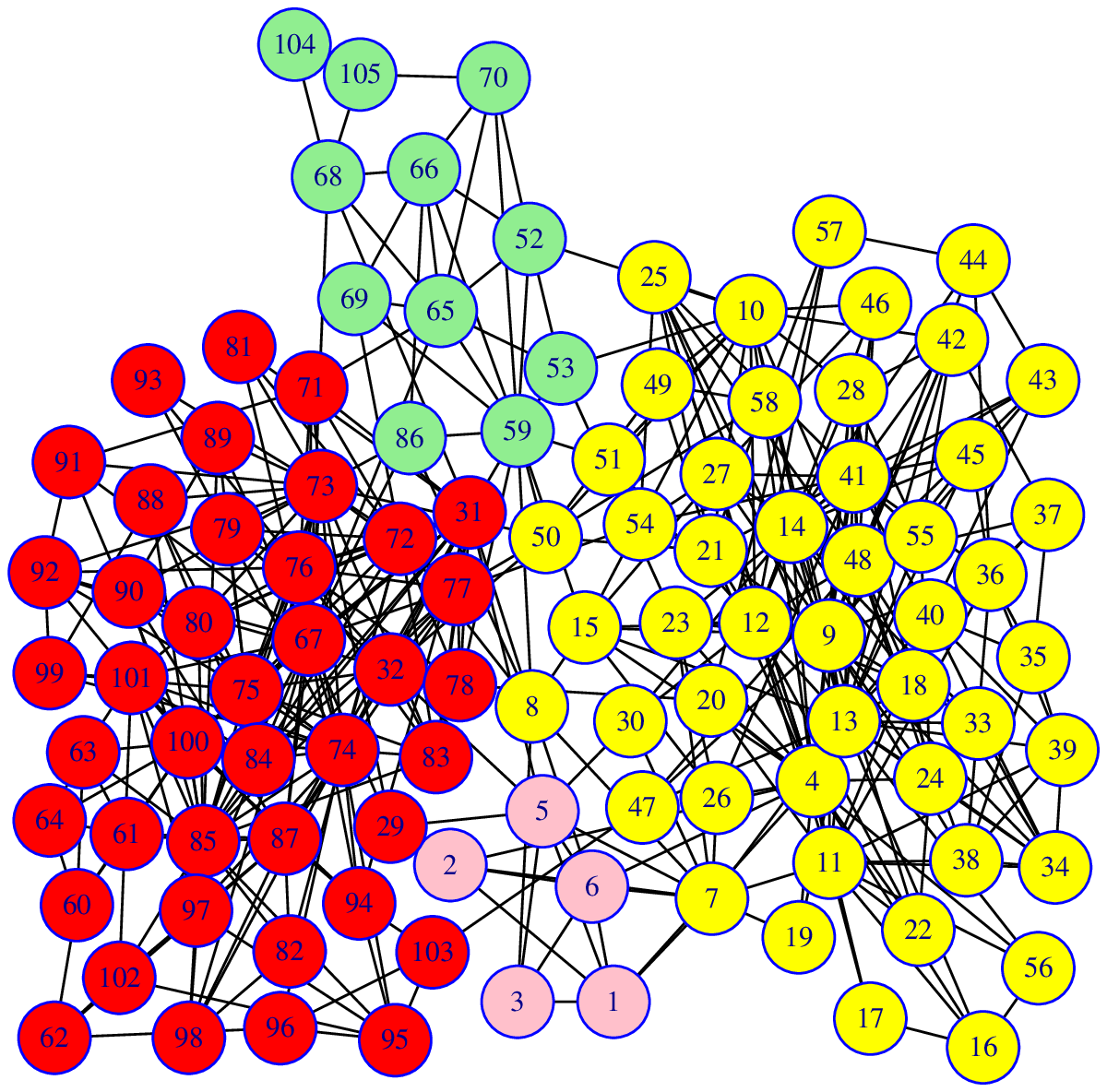}
	    \label{Book_walktrap} \hspace{1pt}%
 	    }
 	    \subfloat[Louvain methodology]{\includegraphics[width=0.5\textwidth]{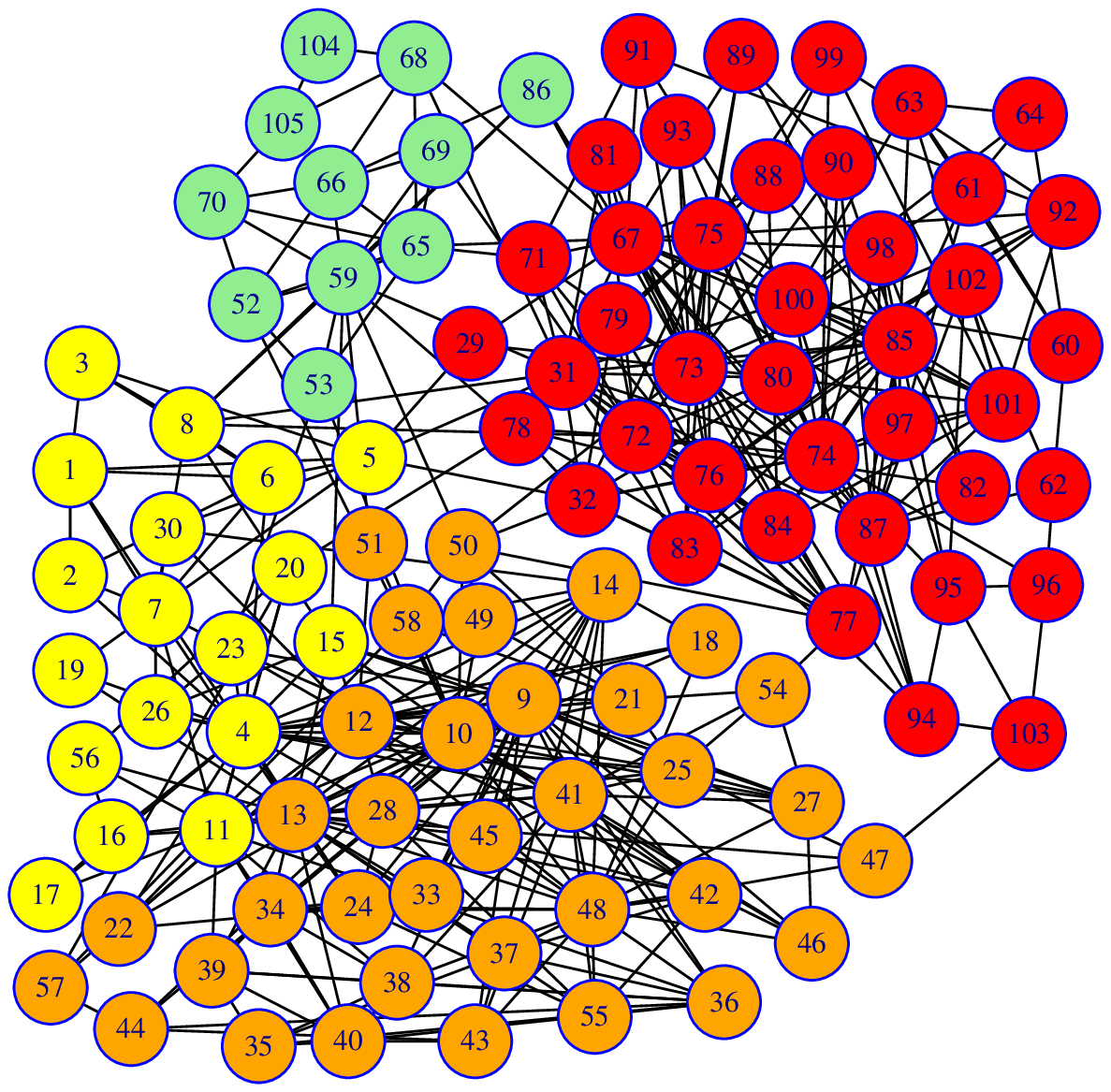}
 	    \label{Books_louvain}
 	    }
	\end{center} 	 
	\caption{Comparison of community detection methods and persistence probability for U.S. political books network.}
\end{figure}

A similar representation is reported in Figure \ref{Book_persistence} for the network of the political books (see Figure \ref{political_books_graph}). We can observe three main peaks of persistence probability, of values $\alpha(V_C)= 0.43$, $\alpha(V_C)= 0.58$ and $\alpha(V_C)= 0.89$, respectively, corresponding to $k=5$, $k=11$ and $k=41$.  Two of these communities corresponds to two communities determined by the Walktrap and Louvain methods as well. They are the two communities with $k=11$, coloured in green in Figures  \ref{Book_walktrap} and \ref{Books_louvain}, and $k=41$, coloured in red. For $k=5$, the third community that we detected is different from what has been found by the two methods. The Walktrap method detects one community with the same size, $k=5$, but different from the one we found, while the Louvain detects nothing. Our methods detects the community $V_C = \{ 60, 61, 63, 64, 100 \}$ while the Walktrap detects $V_C = \{ 1, 2, 3, 5, 6 \}$ (coloured in pink), but it is interesting to observe that for $k=6$ our method detects the similar cluster $V_C = \{1, 2, 3, 5, 6, 7\}$. 

It is worth observing that there is a structural difference between the results of our model and those of both the Walktrap and Louvain method. We do not require that all nodes of the networks belongs to a community, while the latter models do: our method is more flexible without loosing the property of analyzing the network as a whole, while the strict partition could be a too restrictive condition, as it is fully realistic that data reveal that some nodes of the graph are very cohesive between them. The rest of the nodes could be too loosely connected to claim that they are forming exclusive cliques, e.g. they do not show any relevant membership. This is exemplified by the small group detected in the Karate Club, that, at a close inspection, reveals few links outside the group, a property that is not shared by the rest of the nodes. Of course, we are not claiming that what was found by the other methods is wrong, but it is exactly what persistence do: it shows what communities are the most separated from the rest of the graph. This property is further supported through the analysis of the Political books. There, three different groups of books have been found: two of them correspond to clusters also found by the two other methods. This suggests that they are the two most customers' segmentation found by the other methods. On the other side, our method selected a third segment, or cluster, that for numerical reasons has not been detected by the other two. Again, and as for the Karate club, this small group could feature structural properties more consistent than what can be obtained through strict node partitioning.  

\section{Conclusions}\label{conclusions}
    
In this work we presented the integer programming formulation of the problem of finding the node subset with maximum persistence probability and developed heuristic algorithms as well. 
Next, we showed how this methodology helps in discovering communities embedded in a real network by comparison of our findings with well-known methods of community detection. 

There are two main difficulties in applying the persistence index. The first is that the optimal solution is hard to find. Actually, this is a problem shared by many other network statistics, but further research about heuristic procedures is worthwhile. The second difficulty is that the persistence index tends to increase with the subset size $k$ and determining the value of $k$ that corresponds to a community can be problematic. We overcome this issue by locating local peaks on the persistence curve, but further research could be devoted to determine other empirical rule to determine the exact value of $k$. Finally, our model is devoted to finding one community, but it can be used as a subroutine for a graph partitioning model.   

\section*{Acknowledgements}
Research by Stefano Benati has been supported by the financial aid of
NetMeetData: Ayudas Fundacion BBVA a equipos de investigacion cientifica 2019.\\
Research by Rosanna Grassi has been financially supported by Fondo di Ateneo Quota Competitiva 2019, University of Milano - Bicocca.

\bibliographystyle{chicago}
    
\bibliography{References}

\end{document}